\documentclass[11pt]{article}

\usepackage{shapepar,bm}
\usepackage[dvips,arrow,matrix,ps,color,line]{xy}
\usepackage{theorem,amsfonts,amssymb,amscd}
\usepackage{geometry}
\usepackage{graphicx}
\usepackage{graphics}
\usepackage{amssymb}
\usepackage{epstopdf}
\usepackage{amssymb, graphics,hyperref}
\hypersetup{
    colorlinks = true,
    linkcolor = blue,
    anchorcolor = red,
   citecolor = blue,
    filecolor = red,
    pagecolor = red,
    urlcolor = blue}

\newcommand{\mathsym}[1]{{}}

\usepackage[usenames]{color}
\definecolor{MyLightMagenta}{cmyk}{0.1,0.8,0,0.1}
\definecolor{MyDarkBlue}{rgb}{0.1,0,0.3}
\makeindex

\def\NN{\mathbb N}

\def\ZZ{\mathbb Z}

\def\bfm{{\bm\mu}}
\def\bfp{{\mathbf p}}

\def\tta{{\bm\alpha}}

\def\CC{\mathbb C}

\def\PP{\mathbb P}
\def\cocoa{{\hbox{\rm C\kern-.13em o\kern-.07em C\kern-.13em o\kern-.15em A}}}

\def\bfep{{\bm\ep}}

\def\im{\mathrm{Im}}

\def\Acal{{\mathcal A}}

\def\TT{{\bf T}}

\def\Ical{\mathcal I}

\def\ovD{\overline{D}}

\def\ep{{\epsilon}}

\def\wt{{\rm wt}}
\def\ev{{\rm ev}}

\def\w2M{\bigwedge^2M}

\def\wM{\bigwedge M}

\def\w{\wedge }
\def\bw{\bigwedge }

\protect\def\wkM{{\bigwedge^kM}}
\protect
\protect\def\wMn{{\bigwedge M_n}}
\protect

\protect\def\wkMn{{\bigwedge^kM_n}}
\protect

\def\sra{\rightarrow}
\def\lra{\longrightarrow}

\def\proof{\noindent{\bf Proof.}\,\,}
\def\qed{{\hfill\vrule height4pt width4pt depth0pt}\medskip}
\def\be{\begin{equation}}
\def\ee{\end{equation}}
\def\bclm{\begin{claim}}
\def\eclm{\end{claim}}
\def\beqn{\begin{eqnarray}}
\def\eeqn{\end{eqnarray}}
\def\beqn*{\begin{eqnarray*}}
\def\eeqn*{\end{eqnarray*}}

\def\kformep{{\epsilon^1\wedge\ldots\wedge\epsilon^k}}

\def\ikformep{\epsilon^{i_1}\wedge\ldots\wedge\epsilon^{i_k}}

\def\itreformep{\ep^{i_1}\w\ep^{i_2}\w\ep^{i_3}}
\def\iquaformep{\ep^{i_1}\w\ep^{i_2}\w\ep^{i_3}\w\ep^{i_4}}

\def\ihformep{\epsilon^{i_1}\wedge\ldots\wedge\epsilon^{i_h}}

\theoremstyle{change}
\theorembodyfont{\rmfamily}

\newtheorem{claim}{}[section]

\global
\global

\def\no@breaks#1{{\def\\{ \ignorespaces}#1}}    % disallow explicit line breaks

  \makeatletter
\def\cleardoublepage{\clearpage\if@twoside \ifodd\c@page\else
\hbox{} \thispagestyle{empty}
\newpage
\if@twocolumn\hbox{}\newpage\fi\fi\fi} \makeatother

\usepackage{eso-pic,graphicx}
\makeatletter
\newcommand\BackgroundPicture[2]{%
  \setlength{\unitlength}{1pt}%
  default \put(0,\strip@pt\paperheight){%
  \parbox[t][\paperheight]{\paperwidth}{%
    \vfill
     \centering \includegraphics[angle=#2, width=15cm, height=15cm,  bb=0 0 150 150]{#1}
    \vfill
}}} %
\makeatother

\title{Newton Binomial Formulas in Schubert Calculus\footnotemark\footnotetext{{\em Key words and phrases}: Schubert Calculus on a Grassmann algebra, Newton's binomial formulas in Schubert calculus, enumerative geometry of linear series on the projective line; {\em 2000 MSC:} 14M15, 14N15, 15A75.}}

\author{Jorge Cordovez, Letterio Gatto, Ta\'\i se Santiago\thanks{\noindent Work partially sponsored by PRIN ``Geometria sulle Variet\`a Algebriche" (Coordinatore
A.~Verra), INDAM -- GNSAGA, Politecnico di Torino, { FAPESB
(proc. n 8057/2006), CNPq (proc. n 350259/2006-2), UEFS -
Brazil.}}
}
\date{}                                           % Activate to display a given date or no date

\begin{document}
%\thispagestyle{empty}
%\begin{center}
%\includegraphics[height=18.65cm]{cover.pdf}
%\end{center}
%\hfill(To appear on {\em Revista Matem\'atica Complutense})
\maketitle
\hfill{\em {\small To Israel Vainsencher on the occasion of his 60th birthday}\,\,\,\,\,\,\,\,\,\,\,\,\,\,}

\abstract{\small \noindent We prove  Newton's binomial formulas for Schubert Calculus to determine numbers of base point free linear series on the projective line with prescribed ramification divisor supported at given distinct points.}

\section{Introduction}

Let $G(k,n)$ be the complex grassmannian variety parametrizing $k$-dimensional subspaces of $\CC^n$. In~\cite{Gat1}  (see also~\cite{Gat2} and~\cite{Sant1}),  the intersection theory on $G(k,n)$  (Schubert calculus) is   rephrased  via a  natural derivation  on the exterior algebra of a free $\ZZ$-module of rank $n$.   Classical {\em Pieri's} and {\em Giambelli's} formulas are recovered, respectively, from {\em Leibniz's rule} and {\em integration by parts} inherited from such a derivation. The generalization of~\cite{Gat1}  to the   intersection theory on Grassmann bundles is achieved in~\cite{GatSant}, by suitably translating previous important work by Laksov and Thorup (\cite{LakTh}, \cite{LakTh1})  regarding the existence of a (unique)  canonical {\em symmetric structure}  on the exterior algebra of a polynomial ring.

It is natural to ask if the aforementioned  derivation formalism   for Schubert calculus is worthy or if it is nothing more than a mere translation of an old theory into a more or less new language. Indeed,  a couple of years ago, K. Ranestad asked us to test our methods to compute (and possibly to find a formula for) the total number, with multiplicities,  of non projectively equivalent  rational space curves of degree $n+3$ having flexes at given $2n$ distinct  points. Any such a curve is  the image of the morphism $\phi:\PP^1\sra\PP^3$ induced by a very ample $g^3_{n+3}$ on $\PP^1$,  such that the ramification sequence (in the sense of~\cite{Sottile}, Section~1.2) at each of its ramification point   is $(1,2,4,5)$.  Results by Eisenbud and Harris \cite{EHcuspidal} ensure that such a number is finite and equal to the degree of  a suitable product of Schubert cycles -- see Section~\ref{countingwebs} for details. To compute it, we rely on the two main results of this paper, Theorem~\ref{mainnf} and Theorem~\ref{mainnfovd}, regarding certain {\em Newton's  binomial formulas},  which cannot be expressed in the classical  Schubert calculus formalism (as, e.g., in~\cite{GH}). Such formulas allow  us to reduce, after finitely many steps, the degree of {\em any} top codimensional product of  Schubert cycles into an explicit integral linear combination of degrees of Schubert varieties (the latter already  computed by Schubert in~\cite{Sch}).

The Referee, however,  suggested us  that  the same invoked results of~\cite{EHcuspidal} lead naturally to generalize Ranestad's question (and its answer) to that of {\em finding the number, with multiplicities, of $g^r_d$s on $\PP^1$ having pre-assigned ramification at prescribed distinct points}. This last  question was also raised by  Osserman  in~\cite{Osserman}, who conjectures that the scheme of such linear series is reduced when the points are chosen in general position.  In the case $r=1$,   I.~Scherbak proves, via arguments borrowed from representation theory, a nice formula counting  pencils on $\PP^1$ with prescribed ramifications  \cite{Scherb1} (see also~\cite{SV}), while in~\cite{Osserman}, some recursive formulas regarding the case of  higher dimensional linear series  are proven.  Similar counting questions have also been studied by Kharlamov and Sottile in~\cite{Sottile},
 in the more delicate  context of real enumerative algebraic geometry.
 
Our paper   is organized as follows. Section~\ref{nfisc} states, proves and discusses {\em Newton formulas} in Schubert calculus. They remarkably imply that to determine the product structure of $A^*(G(k,n))$ it  is sufficient to know only one Pieri's formula, namely that  expressing  the product of a generator of $A^1(G(k,n))$ with any  Schubert cycle (Cf. remark in Section~\ref{onlys1}). In fact,  the product of any   {\em special Schubert cycle} (see~\cite{GH}) with an arbitrary one, can be  reduced to that primitive case. We do think that this observation is rather new.

In { Section}~\ref{SchubCalc}, Newton's type formulas are then
applied to compute certain powers of {\em special Schubert cycles}
capped with arbitrary classes (modulo rational equivalence) of
Schubert varieties. Such  entirely formal and straightforward
computations  will  finally be exploited  in
Section~\ref{counting} to count  base point free $g^r_d$s on
$\PP^1$, with certain prescribed ramification divisors, for $1\leq
r\leq 3$. In particular,    we use the {\em inclusion-exclusion
formula} (as done in~\cite{SV}   in the context of representation theory)  to supply an alternative proof of (a slightly modified
version of) Scherbak's formula within our framework
(Theorem~\ref{thmscherb}). The method we use seems to be easily
generalizable to the case of grassmannians of higher dimensional
subspaces. However, we shall not follow that path, as we think it is
easier to express an integral  in $A_*(G(k,n))$ as a linear
combination of degrees of Schubert varieties.

Combining the theory exposed in~\cite{EHcuspidal}, or~\cite[Section 1]{Sottile},  with our {\em Newton's formulas} (Theorems~\ref{mainnf} and~\ref{mainnfovd}), we are finally able  to give explicit expressions for the total number $N_{a,b,c,d}$, with multiplicities, of plane irreducible rational curves of degree $n+2$ having $a$ {\em flexes}, $b$ {\em hyperflexes}, $c$ {\em cusps} and $d$ {\em tacnodes} at $a+b+c+d$ distinct points, such that $a+2b+2c+3d=3n$. Similarly, we offer an expression for  the number $f_{a,b,c,d}$ of space rational curves having $a$ {\em stalls}, $b$ {\em hyperstalls} (see~\ref{countingwebs}), $c$ {\em flexes}  and $d$ {\em cusps} at $a+b+c+d$  distinct points  such that $a+2b+2c+3d=4n$. 
%All such expressions  are integral linear combinations of  degrees of Schubert varieties with explicit multinomial coefficients. 
For instance, putting $a=b=d=0$, one gets
\be
f_{0,0,2n,0}=\sum_{\tiny\matrix{c_1+c_2+c_3+c_4=2n\cr 0\leq m\leq c_1+c_2}} {(2n)!\over c_1!c_2!c_3!c_4!}{c_1+c_2\choose m}\omega_{I(c_1,\ldots,c_4;m)},\label{eq:ourranfor}
\ee
where $\omega_{I(c_1,\ldots,c_4;m)}$ is, up to a sign, the degree of a certain Schubert variety explicitly described in Section~\ref{finalsect}. Formula~(\ref{eq:ourranfor}) is  our answer to the original Ranestad's question. In the same vein, putting $a=c=d=0$, one obtains:
\be
HS_n=f_{0,2n,0,0}=\sum_{\tiny\matrix{b_1+\ldots+b_5=2n\cr 0\leq l\leq b_2+b_3}} {(2n)!\over b_1!\cdot \ldots\cdot b_5!}{b_2+b_3\choose l}\omega_{I(b_1,\ldots,b_5;l)},\label{eq:ourrHSfor}
\ee
which is the formula expressing the number (with multiplicities) of (non projectively equivalent) space rational curves having {\em hyperstalls} at prescribed $2n$ distinct points. Indeed, the r.h.s of~(\ref{eq:ourrHSfor}) is  equal to the degree of the $0$-dimensional cycle $\sigma_2^{2n}\cap [G(4,n+4)]$,
 where $\sigma_2$ is the second Chern class of the universal quotient bundle sitting over $G(4,n+4)$.
  One may wish to compute the list of  $HS_n$, $n\geq 0$, by directly evaluating the above degree
  via a mere iteration of Pieri's formula (as, e.g.,  in~\cite{GH}). However,  computations  get very messy already for small values of $n$. Moreover {\em ``Schubert''},  the celebrated {\tt Maple}$^\copyright$ package  designed by Katz and Str{\o}mme \cite{katzstromme},   is apparently unable to go beyond $n=12$  \cite{OkVains}.  Using {\tt Schubert2} on Macaulay2, by Grayson and Stillman \cite{graistil}, one can do much better:  on  a computer  with cpu speed 2.2 Ghz, 16Gb RAM, 4Gb swap,  Jan-Magnus {\O}kland obtained the list of $HS_n$ up to $n=40$, with $HS_{40}$ running for about  8 hours.
 Our formula~(\ref{eq:ourrHSfor}), by contrast, requires no special computer package to be evaluated,  since it is nothing  else than a sum of products of multinomials. In fact, we have been able  to write a trivial \cocoa\, code, successively improved by {\O}kland  \cite{cocoacode},  to write a list of  $HS_n$ via~(\ref{eq:ourrHSfor}). It computes $HS_{42}$ to be
\[
201517182255943002813954873119143476157329393137457696988123090973997900
\]
in about half a hour, showing that the high  computational efficiency of our formulas is  out of reach by the  up to now available packages for doing  Schubert calculus.

 \medskip
{\bf Acknowledgment.} Our first debt of gratitude is with
K.~Ranestad who asked us the main question giving rise to this
paper, and with D.~Laksov for discussions and { constant}
support. For  technical assistance we thank G.~Ghib\`o,
S.~Berrone, A.~Bigatti, but especially  J.~M.~{\O}kland who taught
us the use of Macaulay2 and improved our first  primitive \cocoa\, 
code. We also thank Paolo Aluffi, for useful comments. Finally, we are extremely grateful to the Referee who, with
his/her remarks and suggestions, is responsible for the many substantial
improvements with respect to a previous version of this paper.

\section{Preliminaries}\label{prel}
\claim{} \label{preli1} Let  $X$ be an indeterminate over $\ZZ$,
$M_n:=X\ZZ[X]/(X^{n+1})$ and $\ep^i:=X^i+(X^{n+1})$.  As $m>n$
implies $\ep^m=0_{M_n}$, the $\ZZ$-module  $M_n$ is  freely
generated by  $\bfep^n=(\ep^1,\ldots,\ep^n)$. Let
${\mathbf n}$ be the
set  of the { first $n$ positive integers}.  The {\em weight}
of $I:=(i_1,\ldots,i_k)\in{\mathbf n}^k$ is $\wt(I):=\sum_{j=1}^k(i_j-j)=\sum_{j=1}^k i_j-{1\over
2}k(k+1) $ -- a  non negative integer if the entries of $I$
are all distincts. In $\wkMn$, the $k^{th}$-exterior power of
$M_n$,  we shall often write  $\bfep^I$ { instead } of the
longer expression $\ikformep$.  Let
$I^{\tau}:=(i_{\tau(1)},\ldots,i_{\tau(k)})$, where $\tau\in S_k$,
the symmetric group on $k$ elements: Then
$\bfep^I=\bfep^{I^\tau}$. In particular,
$\bfep^I=0$ if at least two entries of $I$ are equal. If
$\bfep^I\neq 0$, the equality $ \wt(\bfep^I):=\wt(I) $
defines the {\em weight} of $\bfep^I\in\wkMn$.

\claim{} \label{secpkn} Let
$
\Ical^k_n=\{I:=(i_1,\ldots,i_k)\in {\mathbf n}^k\,|\, 1\leq i_1<\cdots<i_k\leq n\}
$
be  the lexicographically (totally) ordered set
of  all the  strictly increasing sequences  of $k$ positive integers not bigger than $n$. We also write $\Ical^{k,w}_n$ for the set of all $I\in\Ical^k_n$ such that $\wt(I)=w$.
Let $\bw^k\bfep^n:=(\bfep^I:I\in\Ical^k_n)$ be the natural  $\ZZ$-basis of $\wkMn$ induced by $\bfep^n$. Denote by  $(\wkMn)_w$  the free submodule of $\wkMn$ generated by the elements of $\bw^k\bfep^n$ of weight $w$. Clearly $(\wkMn)_w=0$ if $w>k(n-k)$. Then $\wkMn$ gets a structure of graded $\ZZ$-module,
$
\wkMn=\bigoplus_{w\geq 0}(\wkMn)_w.
$
\claim{}\label{fundpoint}
The {\em fundamental element} of $\wkMn$ is  the  unique element $g_{k,n}\in\bw^k\bfep^n$ of weight $0$, i.e:
$
g_{k,n}:=\kformep.
$
The {\em point element} of $\wkMn$ is the unique element $\pi_{k,n}\in\bw^k\bfep^n$ of weight $k(n-k)$,  i.e.:
$
\pi_{k,n}:=\ep^{n-k+1}\w\ldots\w\ep^n.
$

\claim{} By \cite{GatSant}, there is a unique sequence  $D:=(D_0,D_1,\ldots)$ of endomorphisms of $\wMn$, the exterior algebra of $M_n$, such that, for each ${\bf p},{\bf q}\in\wMn$,  each $h,l\geq 0$ and each  $i\geq 1$:
\be
\left\{\matrix{ {\displaystyle D_h({\bf p}\w{\bf q})}&=&{\displaystyle \sum_{\tiny{\{h_1,h_2\geq 0\,|\,h_1+h_2=h\}}}D_{h_1}{\bf p}\w D_{h_2}{\bf q}\,\,\,\,\,\,\,\, (h^{th}-Leibniz's \,\,rule);}\cr{}\cr{}\cr\hskip 20pt D_l\ep^i&=&\ep^{i+l}\,\,\hskip 65pt (initial\,\,conditions\,\,on\,\,\bw^1M_n:=M_n).}\right.\label{eq:hleibrule}
\ee
%enjoying  (the $h^{th}$-order) Leibniz's rule:
%\be
%D_h(\alpha\w\beta)=\sum_{\scriptsize{\{h_1,h_2\geq 0\,|\,h_1+h_2=h\}}}D_{h_1}\alpha\w D_{h_2}\beta,\label{eq:hleibrule}
%\ee
%for each $\alpha,\beta\in\wMn$,
%%and the initial conditions $D_h\ep^i=\ep^{i+h}$
% each $h\geq 0$ and each  $i\geq 1$.

Equations~(\ref{eq:hleibrule})  imply that $ D_h(\wkMn)_w\subseteq
(\wkMn)_{w+h} $. Thus $D_h$  is a homogeneous endomorphism of
degree  $0$  of the exterior algebra. The induced  endomorphism of
$\wkMn$, for each $k\geq 1$,  is  homogeneous of degree $h$, with
respect to the weight graduation of $\wkMn$.   The { initial}
conditions~(\ref{eq:hleibrule}) and an easy induction  show that
the endomorphisms  $\{D_i\}_{i\geq 0}$ are pairwise  commuting:
$D_i\circ D_j=D_j\circ D_i$ in $\wMn$. \claim{} For each  $I\in
{\mathbf n}^k$, the {\em degree} of $\bfep^I$,
\[
\int_n \bfep^I=\int_n\ikformep,
\]
 is defined to be $1$ (resp. $-1$)  if there is an even (resp. odd) permutation $\tau\in S_k$ such that $\bfep^{I^\tau}=\pi_{k,n}$, and  zero otherwise. In particular $\int_n\ikformep\neq 0$ implies $\wt(I)=k(n-k)$. By linear extension one gets a function {\em degree} $\int_n: \wkMn\sra\ZZ$, whose kernel is precisely the submodule of $\wkMn$ of elements of weight $< k(n-k)$.

\claim{}\label{gradzt} Let $\ZZ[\TT]:=\ZZ[T_1,T_2,\ldots]$ be the polynomial ring in infinitely many indeterminates. Then $\ZZ[\TT]=\bigoplus_{h\geq 0}\ZZ[\TT]_h$, where one agrees that the degree of a monomial $aT_{i_1}^{j_1}\cdot\ldots\cdot\ldots T_{i_h}^{j_h}$ ($a\in\ZZ$) is  $i_1j_1+\cdots+i_hj_h$. If $I\in\Ical^k_n$,
$
\Delta_I(\TT):=\det(T_{j_i-i})\in \ZZ[\TT]
$
 is the $(I,\TT)$-Schur determinant. It is clearly homogeneous of degree $\wt(I)$.

\claim{}\label{recallsc} The Chow group $A_*(G(k,n))$  of the
grassmannian variety parametrizing $k$-planes in $\CC^n$ is a free
$\ZZ$-module of rank ${n\choose k}$ generated by $\{\Omega_I:
I\in\Ical^k_n\}$, the cycle classes modulo rational equivalence of
Schubert varieties $\Omega_I(E^\bullet)$ associated to some
complete flag $E^\bullet$ of $\CC^n$. The Chow ring $A^*(G(k,n))$
is a $\ZZ$-algebra generated by $\sigma:=(\sigma_i)_{i\geq 0}$,
where $\sigma_i=c_i({\cal Q}_k)$, the $i^{th}$ Chern class of the
(rank $n-k$) universal quotient bundle on $G(k,n)$ (clearly
$\sigma_m=0$ if $m>n-k$). Moreover $A_*(G(k,n))$ is a free module
of rank $1$ over $A^*(G(k,n))$ via {\em cap} product (Poincar\'e
duality). The classical Giambelli's formula can be phrased by saying
that $\Omega_I:=\Delta_I(\sigma)\cap [G(k,n)]$. In particular one
sees that $A^*(G(k,n))$ is a free $\ZZ$-module generated by
$\sigma_I:=\Delta_I(\sigma):={ \det(\sigma_{i_j-i})}$.

\claim{} \label{clm28} Define $D_t:=\sum_{i\geq
0}D_it^i:\wMn\sra\wMn[[t]]$. The first equation of
formula~(\ref{eq:hleibrule}) can be equivalently written as
$D_t({\bf p}\w{\bf q})=D_t{\bf p}\w D_t{\bf q}$ \cite{Gat1}. Let
$D_t^{-1}=\sum_{i\geq 0}(-1)^i\ovD_it^i$ be the formal inverse of
$D_t$ seen as an element of $End_A(\wMn)[[t]]$. Then
\cite{GatSant}: \be D_t^{-1}({\bf p}\w{\bf q})=D^{-1}_t{\bf p}\w
D^{-1}_t{\bf q},\label{eq:hsdrinv} \ee for each ${\bf p},{\bf
q}\in \wMn$.  A direct check shows that~(\ref{eq:hsdrinv})
implies:
\[
\ovD_h({\bf p}\w {\bf q})= \sum_{\tiny{\{h_1,h_2\geq 0\,|\,h_1+h_2=h\}}}\ovD_{h_1}{\bf p}\w \ovD_{h_2}{\bf q},
\]
for each $h\geq 0$ and each ${\bf p},{\bf q}\in \wM$.
For future purposes we observe that the equality $D_t\circ D_t^{(-1)}=1$, holding in $End_A(\wMn)[[t]]$, implies that $\quad \ovD_h=\Delta_{(2\ldots h+1)}(D)$. In particular $\ovD_1=D_1$.
Recall that $h>k$ implies $\wkM\subseteq\ker \ovD_h$ \cite[Proposition 4.1]{GatSant}.

%\claim{\bf Remark.} \label{rmkpartitions} The Schubert cycles $\{\sigma_I, I\in\Ical^k_n\}$ are more commonly indexed by {\em partitions of lenght $\leq k$},  i.e by non increasing sequences of non-negative integers with at most  $k$ non zero terms. More precisely one writes $\sigma_\lambda$, instead of $\sigma_I$,  where $\lambda=(\lambda_1,\ldots,\lambda_k)$ and $\lambda_j:=i_{k+1-j}-(k+1-j)$, for each $1\leq j\leq k$. One usually lists only the non zero terms of a partition $\lambda$ of lenght $\leq k$. So, for instance, within this notation one may write:
%\be
%\ovD_2:=\Delta_{(11)}(D)\qquad\mathrm{and}\qquad \ovD_3=\Delta_{(111)}(D)\label{eq:d11d111}
%\ee
%where $(11)=(3-2,2-1)$ and $(111)=(4-3,3-2,2-1)$.
\claim{} \label{27} Let $E^\bullet$ be any complete flag of
$\CC^n$. If $I\in\Ical^k_n$ and $\wt(I)=w$, \be
\omega_I:=\int_{G(k,n)}\sigma_1^{k(n-k)-w}\cap
\Omega_I=\int_{G(k,n)}\sigma_1^{k(n-k)-w}\sigma_I\cap
[G(k,n)]\label{eq:dshv} \ee is the {\em degree of the Schubert
variety} $\Omega_I(E^\bullet)$ in the Pl\"ucker embedding of
$G(k,n)$. If $J=I^{\tau}$ one defines
\[
\omega_J=sgn(\tau)\cdot \omega_I.
\]
It is known since Schubert \cite{Sch}  that: \be
\omega_I={(k(n-k)-w)!\prod_{j<k} (i_j-i_k)\over
(n-i_1)!\cdot\ldots \cdot (n-i_k)!}.\label{eq:degschvar} \ee
\claim{} \label{dictionary} We denote by $\Acal^*(\wMn)$  the
commutative sub-algebra of $End_\ZZ(\wMn)$, image of the natural
evaluation morphism  $ \ev_D:=\ZZ[\TT]\sra End_\ZZ(\wMn)$, sending
$T_i\mapsto D_i$. There is an obvious restriction morphism
$\rho_k:\Acal^*(\wMn)\sra End_\ZZ(\wkMn)$  mapping $P(D)\in
\Acal^*(\wMn)$ to $P(D)_{|_{\wkMn}}$. Let
$\Acal^*(\wkMn):=\im(\rho_k)$. By~\cite{GatSant}, $\wkMn$ has a
natural structure of  free module of rank $1$ over ${\mathcal
A}^*(\wkMn)$, generated by $g_{k,n}$ (see~\ref{fundpoint}), as a
consequence of Giambelli's formula: \be
\ikformep=\Delta_{(i_1,\ldots,i_k)}(D)\kformep\label{eq:giamb} \ee

Let $\widehat{\Omega}_k:{\wkMn}\sra A_*(G(k,n))$ be the
obvious module isomorphism sending $\bfep^I \mapsto \Omega_I$.
In particular  $\widehat{\Omega}_k(g_{k,n})=[G(k,n)]$, the
fundamental class of $G(k,n)$, and
$\widehat{\Omega}_k(\pi_{k,n})=\Omega_{n-k+1,\ldots,n}=[pt]$, the
class of a point. The main result of~\cite{Gat1}  implies that
there is  a ring isomorphism $\iota_k:A^*(G(k,n))\sra {\mathcal
A}^*(\wkMn)$, sending $\sigma_i\mapsto D_i$, such that  the
diagram:
\begin{small}
\be
\matrix{A^*(G(k,n))\otimes A_*(G(k,n))&\stackrel{\cap}{\lra} &A_*(G(k,n))\cr
{}
\cr
\iota_k\otimes\widehat{\Omega}_k^{-1}\Big\downarrow&{}&\widehat{\Omega}_k\Big\uparrow\cr{}\cr
{\cal A}^*(\bw^kM_n)\otimes \bw^kM_n&{\lra}&\bw^kM_n}\label{eq:diag1}
\ee
\end{small}

\noindent
commutes,
where the bottom horizontal map is defined by the module structure of $\wkMn$ over $\Acal^*(\wkMn)$. In particular,  if  $P\in \ZZ[\TT]$ and
$
\sum_{J\in\Ical^k_n}a_J\cdot \bfep^J
$
is the expansion of $P(D)\bfep^I$ as an integral linear combination of the $\bfep^I$s, then
\[
P(\sigma)\cap \Omega_I=\widehat{\Omega}_k(P(D)\bfep^I)=\widehat{\Omega}_k(\sum_{J\in\Ical^k_n}a_J\bfep^J=\sum_{J\in\Ical^k_n}a_J\widehat{\Omega}_k(\bfep^J)=\sum_{J\in\Ical^k_n}a_J\Omega_J,
\]
where $P(\sigma)\in A^*(G(k,n))$ is the evaluation of $P$ at $\sigma$, via the map $T_i\mapsto \sigma_i$.  Therefore:
\be
\int_{G(k,n)}P(\sigma)\cap \Omega_I=\int_n P(D)\cdot \bfep^I,\label{eq:integ}
\ee
where $\int_{G(k,n)} P(\sigma)\cap \Omega_I$ denotes the  usual degree of the cycle $P(\sigma)\cap \Omega_I\in A_*(G(k,n))$, i.e. the coefficient of $[pt]\in A_*(G(k,n))$, the class of the point of $G(k,n)$, in the expansion of $P(\sigma)\cap \Omega_I$.

\section{Newton's formulas in Schubert Calculus}\label{nfisc}

\claim{} \label{convbin} The binomial coefficient ${n\choose h}$
is, by definition, the coefficient of $a^hb^{n-h}$ in the
expansion of $(a+b)^n$. Therefore ${n\choose j}=0$ if $j<0$, $n<0$
or $j>n$.  If $m, h\geq 0$ are integers, let $
p_h(m)=\{\bfm:=(\mu_1,\ldots, \mu_h)\in \NN^h\,|\, \sum_{i=1}^hm_i=m\} $.
Denote by $|A|$ the cardinality of a set $A$ (we use boldface for multi-indices denoted by greek letters). Then 
\cite[p.~33]{Cameron}:
 \be 
 |p_h(m)|={m+h-1\choose
h-1}.\label{eq:ines} 
\ee 
Define {\em multinomial coefficients} via the equality:
\[
(a_1+\ldots+ a_h)^m=\sum_{\bfm\in p_h(m)}{{m \choose\bfm}}a_1^{\mu_1}\cdot\ldots\cdot a_h^{\mu_h} .
\]
With the usual convention $0!=1$, the  multinomial coefficient can be computed as
\be 
{m \choose \bfm}:={m!\over
\mu_1!\cdot\ldots\cdot\mu_h!} \label{eq:multnm}
\ee
 if $\bfm\in p_h(m)$, while is evidently zero otherwise.
 Equation~(\ref{eq:hleibrule})  for $h=1$,  implies:
\bclm{\bf Proposition.} \label{1stnf}{\em For each
${\bf p},{\bf q}\in\wMn$ and each $h\geq 0$, Newton's binomial
formula holds:
\be D_1^m({\bf p}\w{\bf q})=\sum_{j=0}^m{m\choose
j}D_1^{j}{\bf p}\w D_1^{m-j}{\bf q}.\label{eq:eq1stnf}
\ee }
\eclm
\proof An obvious induction left to the Reader.\qed

 \bclm{\bf Corollary.}\label{cor1stnf} {\em Let $m\geq 0$ and ${\bf p}_1,\ldots,{\bf p}_h\in\wMn$.  Then:
\be
D_1^m({\bf p}_1\wedge\ldots\wedge{\bf p}_h)=\sum_{\bfm\in p_h(m)}{m\choose  \bfm}
D_1^{\mu_1}{\bf p}_1\wedge\ldots\wedge
D_1^{\mu_h}{\bf p}_h.\label{eq:multinew}
\ee
}
\eclm
\proof By induction on
the integer $h\geq 2$. The case $h=2$ is
Proposition~\ref{1stnf}. Suppose  that the formula holds for $h-1$.
Then, by Proposition~\ref{1stnf}:
\begin{eqnarray}
D_1^m({\bf p}_1\wedge{\bf p}_2\wedge\ldots\wedge{\bf p}_h)&=&D_1^m({\bf p}_1\wedge({\bf p}_2\w\ldots\wedge{\bf p}_h))\nonumber\\
&=&\sum_{m_1=0}^m{m\choose m_1}D_1^{m_1}{\bf p}_1\wedge D_1^{m-m_1}({\bf p}_2\wedge\ldots\wedge{\bf p}_h)).\label{eq:lastside}
\end{eqnarray}
By induction, (\ref{eq:lastside}) can be written as:
\begin{eqnarray*}
&{}&\sum_{m_1=0}^m\,\,\sum_{(m_2,\ldots,m_h)\in
p_{h-1}(m-m_1)}{m\choose m_1}{(m-m_1)!\over m_2!\cdot\ldots\cdot m_h!}
D_1^{m_1}{\bf p}_1\w D_1^{m_2}{\bf p}_2\wedge\ldots\wedge D_1^{m_h}{\bf p}_h=\\
 &{}&=\hskip41pt\sum_{\bfm\in p_h(m)}{m\choose \bfm}
 D_1^{\mu_1}{\bf p}_1\wedge\ldots\wedge D_1^{\mu_h}{\bf p}_h.\hskip 41pt \qed
\end{eqnarray*}

\bclm{\bf Lemma.} {\em For each $i\geq 1$, $h\geq 0$ and ${\bf p}\in \wM$, the following formula holds :
\be
D_h(\ep^i\w {\bf p})=\ep^i\w D_h{\bf p}+ D_{h-1}(\ep^{i+1}\w {\bf p}).\label{eq:preip}
\ee
}
\eclm
\proof By a direct check, expanding the two sides of~(\ref{eq:preip}), according to Leibniz's rule~(\ref{eq:hleibrule}).\qed

\bclm{\bf Theorem.}\label{mainnf} {\em Let  $h,m\geq 0$,  $i\geq 1$ and ${\bf p}\in\wMn$. Then:
\be
D_h^{m}( \epsilon^i\w{\bf p}) = \sum ^{m}_{j=0}{m\choose j} D_{h-1}^j( \ep^{i+j}\w D_h^{m-j}{\bf p} )\label{eq:red}
\ee
}
\eclm

\proof  For $m=1$,  formula~(\ref{eq:red}) is~(\ref{eq:preip}). Suppose~(\ref{eq:red}) holds for $m-1$. Since\linebreak
$
D_h^{m}( \epsilon^i\w{\bf p}) = D_h(D_h^{m-1}( \epsilon^i\w{\bf p})),
$
induction on $m$ gives:
\[
D_h^{m}(\epsilon^i\w{\bf p})=D_h\left(\sum ^{m-1}_{j=0}{m-1\choose j}D_{h-1}^j(\epsilon^{i+j}\w D^{m-1-j}_h{\bf p})\right).
\]
Using the linearity and the fact that the operators $\{D_h\}_{h\geq 0}$ are pairwise commuting:
\[
D_h^{m}(\ep^i\w{\bf p})=\left(\sum ^{m-1}_{j=0}{m-1\choose j}D_{h-1}^jD_h(\epsilon^{i+j}\w D^{m-1-j}_h {\bf p})\right),
\]
from which, by  applying~(\ref{eq:preip}) again, one gets:
\begin{small}
\begin{eqnarray*}
D_h^{m}(\epsilon^i\w{\bf p} )=
\sum ^{m-1}_{j=0}{m-1\choose j}D_{h-1}^j(\ep^{i+j}\w D_h^{m-j}{\bf p}+D_{h-1}( \epsilon^{i+j+1}\w D^{m-1-j}_h {\bf p}))=\\
=\sum ^{m-1}_{j=0}{m-1\choose j}D_{h-1}^j(\ep^{i+j}\w D_h^{m-j}{\bf p})+\sum ^{m-1}_{j=0}{m-1\choose j}D_{h-1}^{j+1}(\ep^{i+j+1}\w D_h^{m-1-j}{\bf p})=\\
=\sum ^{m}_{j=0}{m-1\choose j}D_{h-1}^j(\ep^{i+j}\w D_h^{m-j}{\bf p})+\sum ^{m}_{j=0}{m-1\choose j-1}D_{h-1}^j(\ep^{i+j}\w D_h^{m-j}{\bf p})=\\
=\sum ^{m}_{j=0}\left[ {m-1\choose j}+{m-1\choose j-1} \right] D_{h-1}^j(\ep^{i+j}\w D_h^{m-j}{\bf p})=\sum ^{m}_{j=0}{m\choose j} D_{h-1}^j( \ep^{i+j}\w D_h^{m-j}{\bf p} ).
\end{eqnarray*}
\qed
\end{small}

\bclm{\bf Theorem.}\label{mainnfovd} {\em  Let ${\bf p}\in\wMn$. Then for each $i\geq 1$ and each $h\geq 0$ one has (notation as in~\ref{clm28}):
\be
\ovD_h^m(\ep^i\w\bfp)=\sum_{j=0}^m{m\choose j}\ep^{i+j}\w \ovD_{h-1}^j\ovD_h^{m-j}\bfp\label{eq:redovd}
\ee
}
\eclm
\proof
By induction on $m$.  Recall that by~\cite{GatSant}, Proposition~4.1, $\ovD_j\ep^i=0$ unless $0\leq j\leq 1$, in which case one has $\ovD_0\ep^i=\ep^i$ and $\ovD_1\ep^i=\ovD_1\ep^i=\ep^{i+1}$. Then:
\[
\ovD_h(\ep^{i}\w\bfp)=\sum_{j=0}^h\ovD_j{ \ep^{i}\w
\ovD_{h-j}}\bfp=\ep^i\w \ovD_h\bfp+\ep^{i+1}\w\ovD_{h-1}\bfp
\]
i.e.~(\ref{eq:redovd}) holds for $m=1$.
 Suppose it holds for $m-1$. Then
\[
\ovD_h^m(\ep^i\w\bfp)=\ovD_h(\ovD_h^{m-1}(\ep^i\w\bfp))=\ovD_h\sum_{j=0}^{m-1}{m-1\choose j}\ep^{i+j}\w \ovD_{h-1}^{j}\ovD_h^{m-1-j}\bfp=
\]
i.e. by linearity and the definition of $\ovD_h$:
\[
=\sum_{j=0}^{m-1}{m-1\choose j}\ep^{i+j}\w\ovD_{h-1}^{j}\ovD_h^{m-j}\bfp+ \sum_{j=0}^{m-1}{m-1\choose j}\ep^{i+j+1}\w\ovD_{h-1}^{j+1}\ovD_h^{m-1-j}\bfp=
\]
\[
=\sum_{j=0}^{m}{m-1\choose j}\ep^{i+j}\w\ovD_{h-1}^{j}\ovD_h^{m-j}\bfp+ \sum_{j=0}^{m}{m-1\choose j-1}\ep^{i+j}\w\ovD_{h-1}^{j}\ovD_h^{m-1}\bfp=
\]
\[
\sum_{j=0}^{m}\left[{m-1\choose j}-{m-1\choose j-1}\right]\ep^{i+j}\w\ovD_{h-1}^{j}\ovD_h^{m-j}\bfp=\sum_{j=0}^{m}{m\choose j}\ep^{i+j}\w\ovD_{h-1}^{j}\ovD_h^{m-j}\bfp
\]
as claimed.
\qed
\claim{} \label{onlys1} Formula~(\ref{eq:red}) supports an explicit efficient algorithm to express any  product of special Schubert cycles as an  integral linear combination of products of the form   $\sigma_1^m\cap\sigma_I$. Such an algorithm is extremely useful to perform computations, but it has also a theoretical relevance: it shows that the algebra structure of $A^*(G(k,n))$ is completely determined once one knows the product $\sigma_1\sigma_I$ for each $I\in\Ical^k_n$. To this purpose we exploit diagram~(\ref{eq:diag1}) together with our main formula~(\ref{eq:red}), as follows.
For each $I:=(i_1,\ldots,i_k)\in\Ical^k_n$, $h\geq 1$ and $d,m\geq 0$, let
\[
{\cal J}(d,k,I):=\{(j,J)\in \NN\times {\Ical^k_n}\,|\,   j+\wt(J)=\wt(I)+d\}.
\]

\bclm{\bf Proposition.}\label{propred} {\em There is  an explicit
algorithm to express, in at least one way, any
$D^m_h\big(\bfep^I\big)$  as a $\ZZ$-linear combination of the
elements of the set \be \{D_1^j\big(\bfep^J\big)\,|\,
(j,J)\in{\cal J}(mh,k,I)\}\label{eq:setlc} \ee i.e: $
D^m_h\big(\bfep^I\big)=\sum_{(j,J)\in{\cal
J}(mh,k,I)}a_JD_1^j\big(\bfep^J\big). $ } \eclm \proof
Induction on $h\geq 1$ and on $k\geq 2$. For $h=1$ the proposition
is trivial for all $k\geq 2$.  Let us assume that it holds for all
$1\leq h'\leq h-1$, and  all $\ep^{i_1}\w\ep^{i_2}\in \bw^2M_n$
($k=2$).  By~(\ref{eq:red}):
\[
D_h^m(\ep^{i_1}\w\ep^{i_2})=\sum_{j=0}^m{m\choose j}D_{h-1}^j\big(\ep^{i_1+j}\w \ep^{i_2+h(m-j)}\big).
\]
By the inductive hypothesis $D_{h-1}^j(\ep^{i_1+j}\w
\ep^{i_2+h(m-j)})$ is equal to a suitable linear combination of
elements of the form $D_1^{ j'}(\ep^{j_1}\w\ep^{j_2})$. Hence
the proposition holds for $k=2$ and for all $h\geq 1$. Suppose now
that the proposition holds i) for  all $1\leq h'\leq h-1$ and  all
$I\in\Ical^k_n$ ($k\geq 2$) and ii) for all $h\geq 1$,  and all
$I':=(i_2,\ldots,i_k)\in{\Ical^{k-1}_n}$. We apply~(\ref{eq:red})
to $\ep^{i_1}\w {\bf p}$, with  ${\bf
p}:=\bfep^{I'}=\ep^{i_2}\w\ldots\w\ep^{i_k}$: 
\be
D_h^m(\ikformep)=\sum_{j=0}^m{m\choose j}D_{h-1}^j(\ep^{i_1+j}\w
D_h^{m-j}(\ep^{i_2}\w\ldots\w\ep^{i_k})).\label{eq:red3} \ee By
induction,  $ D_h^{m-j}(\ep^{i_2}\w\ldots\w\ep^{i_k})$ can be
written as:
\[
\sum_{j=0}^m\,\,\sum_{(j',J')\in{\cal J}(h(m-j),k-1,I')} a_{J'}D_1^{j'}\big(\bfep^{J'}\big),
\]
for a suitable choice of $a_{J'}\in\ZZ$, and~(\ref{eq:red3}) takes the form:
\[
D_h^m(\ikformep)=\sum_{j=0}^m{m\choose j}D_{h-1}^j(\ep^{i_1+j}\w \sum_{(j',{J'})\in{\cal J}(h(m-j),k-1,I')} a_{J'}D_1^{j'}\big(\bfep^{J'}\big) ).
\]
By Corollary~\ref{cor1stnf} one has:
\[
D_1^{j'}\big(\bfep^{J'}\big)=\sum_{(l_2,\ldots,l_{k})\in p_{k-1}(j')}{j'!\over  l_2!\cdot \ldots\cdot l_{k}!}\ep^{j'_2+l_2}\w\ldots\w\ep^{j'_k+l_k}.
\]
In conclusion:
\[
D_h^m(\ikformep)=\sum_{j=0}^m{m\choose
j}D_{h-1}^j(\sum_{(l_2,\ldots,l_{k}){ \in}
p_{k-1}(j')}{j'!\over l_2!\cdot\ldots\cdot
l_{k}!}\ep^{i_1}\w\ep^{j'_2+l_2}\w\ldots\w\ep^{j'_k+l_k}=
 ).
\]

\[
=\sum_{\tiny\matrix{0\leq j\leq m\cr (l_2,\ldots,l_{k}){ \in}
p_{k-1}(j')}}{m\choose j}{j'!\over  l_2!\cdot \ldots\cdot
l_{k}!}D_{h-1}^j(\ep^{i_1}\w\ep^{j'_2+l_2}\w\ldots\w\ep^{j'_k+l_k}
).
\]
By the inductive hypothesis, one then sees that $D_h^m(\ikformep)$ is itself an integral linear combination of the elements of the set~(\ref{eq:setlc}).\qed
\bclm{\bf Corollary.}\label{corind} {\em Let $P\in\ZZ[\TT]_h$  (see Section~\ref{gradzt}). There is an explicit effective algorithm to express $P(D)\ikformep$, in at least one way,   as   an integral linear combination of
$
\{D_1^j(\bfep^J)\}_{(j,J)\in{\cal J}(h,k,I)}
$.
}
\eclm
\proof
Any such $P$ is, by definition, a (finite) integral linear combination of monomials of the form $T_{h_1}^{m_1}T_{h_2}^{m_2}\cdot\ldots\cdot T_{h_l}^{m_l}$, and then  it suffices  to check the property for any such a term. One argues by induction on the integer $l$: for $l=1$ the property holds. If it holds for $l-1$, one has that
\[
D_{h_1}^{m_1}D_{h_2}^{m_2}\cdot\ldots\cdot D_{h_l}^{m_l}(\ikformep)=D_{h_1}^{m_1}(D_{h_2}^{m_2}\cdot\ldots\cdot D_{h_l}^{m_l}(\ikformep))
\]
By induction, there are integers $a_J$ such that
\[
D_{h_2}^{m_2}\cdot\ldots\cdot D_{h_l}^{m_l}(\ikformep)=\sum_{(j,J)}a_JD_1^j\bfep^J
\]
the sum over all $(j,J)$ such that $j+\wt(J)=wt(I)+h_2m_2+\cdots+h_lm_l$.
Hence:
\[
D_{h_1}^{m_1}D_{h_2}^{m_2}\cdot\ldots\cdot D_{h_l}^{m_l}(\ikformep)=D_{h_1}^{m_1}(\sum_{(j,J)}a_JD_1^j\bfep^J)=\sum_{(j,J)}a_JD_1^j(D_{h_1}^{m_1}\bfep^J)
\]
and one finally concludes by applying  Proposition~\ref{propred} again.\qed
\bclm{\bf Corollary.}\label{cor39} {\em Suppose $\deg(P)=k(n-k)$. There is an explicit algorithm to compute
\[
\int_{G(k,n)}P(\sigma)\cap[G(k,n)]=\int_nP(D)(\kformep)
\]
as a $\ZZ$-linear combination of degrees of Schubert varieties.
}
\eclm
\proof
In fact, by Corollary~(\ref{corind})  one can determine integers $a_J$ to write:
\[
P(D)\kformep=\sum_{(j,J)}a_JD_1^j\bfep^J,
\]
where in each summand $j+\wt(J)=k(n-k)$.
Taking integrals:
\[
\int_nP(D)\kformep=\sum_{(j,J)}a_J\int_nD_1^j\bfep^J.
\]
Since $j+\wt(J)=k(n-k)$,  $\int_nD_1^jJ\bfep^J$ is precisely the degree of the Schubert variety $\Omega_J(E^\bullet)$, $E^\bullet$ being an arbitrary complete flag of $\CC^n$ (Cf. Sections~\ref{recallsc}, \ref{27}).\qed

\section{\bf Computations in $\bw^kM_{n}$, $2\leq k\leq 4$ }\label{SchubCalc}
By~\cite{EHcuspidal}, Schubert calculus on $G(k,n+k)$ can be interpreted in terms of enumerative geometry of linear series on the projective line with prescribed ramification divisor. This fact motivates the computation below which shall be applied to enumerative problems regarding rational curves in next section.

\claim{} Notation as in Section~\ref{clm28}. Recall that for each $h\geq 0$ \cite[Sect.~2.10]{GatSant}:
\[
\ovD_h(\ihformep)=\ep^{i_1+1}\w\ldots\w\ep^{i_h+1},
\]
and hence, by induction:
\be
\ovD_h^m(\ihformep)=\ep^{i_1+m}\w\ldots\w\ep^{i_h+m}.\label{eq:iterdbar}
\ee
\bclm{\bf Proposition.}\label{dbarhm1} For each $h\geq 1$,
\be
\ovD_{h-1}^m(\ihformep)=\sum_{\bfm\in p_h(m)}{m\choose \bfm}\ep^{i_1+\sum_{j\neq 1}\mu_j}\w\ldots\w\ep^{i_h+\sum_{j\neq h}\mu_j}\label{eq:lemcomp}
\ee
\eclm
\proof
If $h=1$, $D_0\ep^{i_1}=\ep^{i_1}$ and the proposition is trivial.
For $h=2$, is just Newton formula~(\ref{eq:eq1stnf}), since $\ovD_1=D_1$ and
\[
D_1^m(\ep^{i_1}\w\ep^{i_2})=\sum_{m_1=0}^m {m\choose m_1}\ep^{i_1+m_1}\w\ep^{i_2+m-m_1}=\sum_{\bfm\in p_2(m)}{m\choose \bfm}\ep^{i_1+\mu_1}\w\ep^{i_2+\mu_2}.
\]
Suppose the formula true for $h-1$. Then  
\[
\ovD_{h-1}^m(\ihformep)=\ovD_{h-1}^m(\ep^{i_1}\w(\ep^{i_2}\w\ldots\w\ep^{i_h}))=
\]
and by~(\ref{eq:redovd}):
\[
=\sum_{m_1=0}^m{m\choose m_1}\ep^{i_1+m_1}\w
\ovD^{m_1}_{h-2}(\ep^{i_2+m-m_1}\w\ldots\w\ep^{i_h+m-m_1})=
\]
\[
=\sum {m\choose m_1}{m_1!\over m_2!\cdot \ldots \cdot
m_{h-1}!}\ep^{i_1+m_1}\w\ep^{i_2+m-m_1{ +}\sum_{j\neq
2}m_j}\w\ldots\w\ep^{i_h+m-m_1+\sum_{j\neq h}m_j}=
\]
\[
=\sum{m!\over (m-m_1)!\cdot m_2!\cdot\ldots\cdot
m_h!}\,\ep^{i_1+m_1}\w\ep^{i_2+m-m_1{ +}\sum_{j\neq
2}m_j}\w\ldots\w\ep^{i_h+m-m_1+\sum_{j\neq h}m_j},
\]
where the last two sums are over all $(m_1,m_2,\ldots, m_h)$ such that $0\leq m_1\leq m$ and $\sum_{j=2}^hm_j=m_1$.
Taking  $\bfm\in p_h(m)$ such that $\mu_1=m-m_1$ and $\mu_j=m_j$, for $2\leq j\leq h$, so that $\mu_1=m_2+\cdots+m_h$, one gets exactly formula~(\ref{eq:lemcomp}). \qed

Recall the notation of Section~\ref{preli1}.
\bclm{\bf Proposition.} {\em Let $n\geq 0$ and $M_{n}$ as in \ref{preli1}  and $D:=(D_1,D_2,\ldots)$ as in formula~(\ref{eq:hleibrule}). Then the following equalities holds in $\wMn$:
 \be
\mathrm{i)}\hskip9pt D_2^m(\itreformep)=\sum_{\bfm\in p_4(m)}{m\choose \bfm}D_1^{\mu_1}(\ep^{i_1+\mu_1}\w\ep^{i_2+\mu_2+2\mu_4}\w\ep^{i_3+2\mu_3+\mu_2});\hskip9pt\label{eq:forD2mw3}
\ee

\be
\mathrm{ii)}\hskip15pt D_2^m(\ep^{i_1}\w\ep^{i_2}\w\ep^{i_3}\w\ep^{i_4})=\sum_{\tiny\matrix{\bfm\in p_5(m)\cr  0\leq l\leq \mu_2)+\mu_3}}{m\choose \bfm}{\mu_2+\mu_3\choose l}D_1^{\mu_1}(\bfep^{I(\bfm,{l})}),\hskip8pt\label{eq:forD2m}
\ee
where
\[
I(\bfm,l)=(i_1+\mu_1,i_2+\mu_2+2\mu_5,i_3+\mu_3+l,i_4+\mu_2+2\mu_4+\mu_3-l);
\]

\be
\mathrm{iii)}\hskip15pt \ovD_2^m(\ep^{i_1}\w\ep^{i_2}\w\ep^{i_3}\w\ep^{i_4})=\sum_{\tiny\matrix{\bfm\in p_4(m)\cr  0\leq l\leq\mu_1+\mu_2}}{m\choose \bfm}{\mu_1+\mu_2\choose l}\bfep^{J(\bfm,{l})},\hskip8pt\label{eq:forovD2m}
\ee
where
\[
J(\bfm,l):=(i_1+\mu_1+\mu_4,i_2+\mu_2+\mu_4, i_3+\mu_3+l, i_4+\mu_1+\mu_2+\mu_3-l).
\]

}
\eclm

\proof
It consists in a repeated application of~(\ref{eq:red}). We limit ourselves to  the verification of~(\ref{eq:forD2m}) and~(\ref{eq:forovD2m}), leaving the easier~(\ref{eq:forD2mw3}) to the reader's care, as a more or less amusing exercise.

\medskip
Regarding formula~(\ref{eq:forD2m}) we first observe that:
\[
D_2^m(\iquaformep)=\sum_{\alpha_1=0}^m{m\choose \alpha_1}D_1^{\alpha_1}(\ep^{i_1+\alpha_1}\w D_2^{m-\alpha_1}(\ep^{i_2}\w\ep^{i_3}\w\ep^{i_4}))=
\]
\[
=\sum_{\tiny\matrix{0\leq a_1\leq m\cr 0\leq a_2\leq
m-\alpha_1}}^m{m\choose \alpha_1}{m-\alpha_1\choose a_2}D_1^{\alpha_1}(\ep^{i_1+\alpha_1}\w
D_1^{\alpha_2}(\ep^{i_2+\alpha_2}\w D_2^{m-\alpha_1-\alpha_2}{ (
\ep^{i_3}\w\ep^{i_4})}))=
\]

\be
=\sum_{{\bm\alpha}\in p_4(m)}{m\choose {\bm\alpha}} D_1^{\alpha_1}(\ep^{i_1+\alpha_1}\w D_1^{\alpha_2}(\ep^{i_2+\alpha_2}\w D_1^{\alpha_3}(\ep^{i_3+\alpha_3}\w \ep^{i_4+2\alpha_4}))\label{eq:initD2m}
\ee
where~(\ref{eq:initD2m}) has been gotten by repeatedly applying~(\ref{eq:red}) and having set ${\bm\alpha}=(\alpha_1,\alpha_2,\alpha_3,\alpha_4)$.
Now, by~(\ref{eq:eq1stnf}):
\[
D_1^{\alpha_2}(\ep^{i_2+\alpha_2}\w D_1^{\alpha_3}(\ep^{i_3+\alpha_3}\w \ep^{i_4+2\alpha_4}))=
\]
\[
=\sum_{0\leq b\leq \alpha_2}{\alpha_2\choose b}\ep^{i_2+\alpha_2+b}\w D_1^{\alpha_3+\alpha_2-b}(\ep^{i_3+\alpha_3}\w \ep^{i_4+2\alpha_4}).
\]
By {definition} of $D_1$ and by applying~(\ref{eq:eq1stnf})
once more:
\[
\sum_{\tiny\matrix{0\leq b\leq \alpha_2\cr 0\leq l\leq \alpha_3+\alpha_2-b}}{\alpha_2\choose b}{\alpha_3+\alpha_2-b\choose l}\cdot\bfep^{J_1({\bm\alpha}, b, l)}
\]
with $J_1(\tta,b,l):=(i_2+\alpha_2+b,i_3+\alpha_3+l,i_4+2\alpha_4+\alpha_2+\alpha_3-b-l)$.

Last expression plugged into~(\ref{eq:initD2m}) gives:
\be
\sum_{\tiny\matrix{\tta\in p_4(m)\cr 0\leq b\leq \alpha_2\cr
0\leq l\leq \alpha_3+\alpha_2-b}}{ {m\choose 
\tta}{{\alpha_2}\choose b}{{\alpha_3+\alpha_2-b}\choose
l}}\cdot D_1^{\alpha_1}\bfep^{J_2(\tta,b,l)},\label{eq:912j2}
\ee
where $J_2(\tta,b,l):=(i_1+\alpha_1,i_2+\alpha_2+b,i_3+\alpha_3+l,i_4+2\alpha_4+\alpha_2+\alpha_3-b-l)$,
i.e., using the definition of the multinomial coefficient:
\[
\sum{\tiny {m!\over {
\alpha_1!(\alpha_2-b)!\alpha_3!\alpha_4!b!}}{{\alpha_3+\alpha_2-b}\choose
l}}\cdot D_1^{\alpha_1}\big(\bfep^{J_2(\tta, b,l)}\big).
\]

We may rename the variables defining $\bfm\in p_5(m)$ as (the way to do that is not unique)
$
\mu_1:=\alpha_1;\quad \mu_2=\alpha_2-b;\quad  \mu_3=\alpha_3; \quad \mu_4=\alpha_4;\quad \mu_5=b
$.
Finally, expressing $\tta,b$ as functions of  such  a $\bfm\in p_5(m)$, and substituting in~(\ref{eq:912j2}), one gets:
\[
\sum_{\tiny\matrix{\bfm\in p_5(m)\cr 0\leq l\leq \mu_2+\mu_3}}{m\choose \bfm}{\mu_2+\mu_3\choose l}D_1^{\mu_1}(\bfep^{I(\bfm,l)}),
\]
having set :
$
I(\bfm,l)=(i_1+\mu_1,i_2+\mu_2+2\mu_5,i_3+\mu_3+l,i_4+\mu_2+2\mu_4+\mu_3-l),
$
which is precisely~(\ref{eq:forD2m}).

We come now to the proof of~(\ref{eq:forovD2m}). By  applying  formula~(\ref{eq:redovd}) once,  one gets:
\[
\ovD_2^m(\iquaformep)=\sum_{\alpha_1=0}^m{m\choose \alpha_1}\ep^{i_1+\alpha_1}\w \ovD_1^{\alpha_1}\ovD_2^{m-\alpha_1}(\ep^{i_2}\w\ep^{i_3}\w\ep^{i_4})=
\]
i.e., again:
\[
=\sum_{\tiny\matrix{0\leq \alpha_1\leq m\cr 0\leq \alpha_2\leq
m-\alpha_1}}{m\choose \alpha_1}{m-\alpha_1\choose \alpha_2}\ep^{i_1+\alpha_1}\w {
\ovD_1}^{\alpha_1}(\ep^{i_2+\alpha_2}\w\ovD_1^{\alpha_2}\ovD_2^{m-\alpha_1-\alpha_2}(\ep^{i_3}\w\ep^{i_4}))=
\]
\be
=\sum_{\tiny\matrix{\alpha_i\geq 0\cr \alpha_1+\alpha_2+\alpha_3=m}}{m!\over
\alpha_1!\alpha_2!\alpha_3!}\ep^{i_1+\alpha_1}\w{
\ovD_1}^{\alpha_1}(\ep^{i_2+\alpha_2}\w\ovD_1^{\alpha_2}(\ep^{i_3+\alpha_3}\w\ep^{i_4+\alpha_3}))=\label{eq:904cns}
\ee
where in the last equality we used~(\ref{eq:multnm}) (having put $\alpha_3:=m-\alpha_1-\alpha_2$) and~(\ref{eq:iterdbar}).
By applying twice
formula~(\ref{eq:eq1stnf}) to  expression~(\ref{eq:904cns}) one gets:
\[
\sum_{\tiny\matrix{\alpha_1+\alpha_2+\alpha_3=m\cr 0\leq b\leq \alpha_1}}{m!\over \alpha_1!\alpha_2!\alpha_3!}{\alpha_1\choose b}\ep^{i_1+\alpha_1}\w \ep^{i_2+\alpha_2+b}\w\ovD_1^{\alpha_2+\alpha_1-b}(\ep^{i_3+\alpha_3}\w\ep^{i_4+\alpha_3})=
\]

\[
=\sum_{\tiny\matrix{\alpha_1+\alpha_2+\alpha_3=m\cr 0\leq b\leq \alpha_1\cr 0\leq l\leq \alpha_2+\alpha_1-b}}{m!\over \alpha_1!\alpha_2!\alpha_3!}{\alpha_1\choose b}{\alpha_2+\alpha_1-b\choose l}\ep^{i_1+\alpha_1}\w \ep^{i_2+\alpha_2+b}\w\ep^{i_3+\alpha_3+ l}\w\ep^{i_4+\alpha_3+\alpha_1+\alpha_2-b- l}=
\]
\[
=\sum_{\tiny\matrix{\alpha_1+\alpha_2+\alpha_3=m\cr 0\leq b\leq \alpha_1\cr 0\leq l\leq \alpha_2+\alpha_1-b}}{m!\over (\alpha_1-b)!\alpha_2!\alpha_3!b!}{\alpha_2+\alpha_1-b\choose  l}\ep^{i_1+\alpha_1}\w \ep^{i_2+\alpha_2+b}\w\ep^{i_3+\alpha_3+ l}\w\ep^{i_4+\alpha_3+\alpha_1+\alpha_2-b- l}.
\]
Putting $\bfm\in p_4(m)$ such that $\mu_1=\alpha_1-b$, $\mu_2=\alpha_2$, $\mu_3=\alpha_3$ and $\mu_4=b$, one finally obtains~(\ref{eq:forovD2m}), as desired.\qed
\section{Counting $g^r_{n+r}$ on $\PP^1$ with prescribed ramification at distinct points.}\label{counting}
The main references for next subsection are ~\cite{EHcuspidal}  and~\cite[Section 1]{Sottile}.
\claim{}
If $a$ is a non-negative integer and  $W$  a subspace of $F:=H^0(\PP^1,O_{\PP^1}(d))$, denote by $W(-aP)$  the subspace of all the sections of $W$ vanishing at $P$ with multiplicity at least $a$. A $g^r_d$ on $\PP^1$ is the choice of a $(r+1)$-dimensional subspace $V$ of $F$.
Since $F(-aP)$ has codimension $a$,  for each $0\leq a \leq d+1$, the chain of inclusions
\[
F^\bullet(P): F\supset F(-P)\supset\cdots\supset F(-dP)\supset F(-(d+1)P)=0,
\]
defines the (complete)  {\em osculating flag} of $F$ at $P$.
 If there exists a non zero $v\in V$ vanishing at $P$ with multiplicity exactly $i-1$ ($i\geq 1$), one says that $i-1$ is a $V$-{\em order} at $P$. There are exactly $r+1$ orders at each point $P$,  forming its $V$-{\em order sequence} $0\leq i_1-1<\ldots<i_{r+1}-1\leq d$. Following~\cite[Section~1.2]{Sottile}, we shall say that $1\leq i_1<\ldots<i_{r+1}\leq d+1$ is the $V$-{\em ramification sequence} at $P$. The $V$-{\em weight} of $P$ is $\wt_V(P)=\sum_{j=1}^{r+1}(i_j(P)-j)$. A point $P\in\PP^1$ is a $V$-ramification point if one of the following equivalent conditions occur: i)  $\dim V(-(r+1)P)>1$, ii)  $\wt_V(P)>0$, iii) the $V$-ramification sequence at $P$ is different from $(1,2,3,\ldots,r+1)$.
All but finitely many points have weight $0$ and the {\em total weight} of the ramification points is prescribed by the Brill-Segre-Pl\"ucker formula:
\[
\wt_V=\sum_{P\in C}\wt_V(P)=(r+1)(d-r).
\]
\claim{} \label{recalEH} If  $((i_1(P),\ldots, i_{r+1}(P))$ is the $V$-ramification sequence at $P$,  then $V$ belongs,  by \cite{EHcuspidal}, to the Schubert variety $\Omega_{i_1,\ldots, i_{r+1}}(F^\bullet(P))$ of the grassmanniann $G(r+1, d+1)$. Moreover, if $P_1,\ldots, P_m$ are $m$ distinct points and $I_1,\ldots,I_{m}\in \Ical^{r+1}_d$, the Schubert varieties $\Omega_{I_j}(F^\bullet(P_j))$ are {\em dimensionally transverse}, i.e. every irreducible component of $\bigcap_{j=1}^m\Omega_{I_j}(F^\bullet(P_j))$ has codimension $\sum{\rm codim}\, \Omega_{I_j}(F^\bullet(P_j))=\sum_{j=1}^m\wt(I_j)$. In particular, if the intersection is zero dimensional, the degree of the intersection cycle is given by the integral of the product of the corresponding Schubert cycles, in the intersection ring of $G(r+1, d+1)$, capped with the fundamental class of $G(r+1, d+1)$. The formulas gotten in Section~\ref{SchubCalc}  via~(\ref{eq:red}) can be easily applied to get expressions for the number (with multiplicity) of  $g^r_{n+r}$ on $\PP^1$  with prescribed ramifications at prescribed points  for $r=1,2,3$, i.e. for {\em pencils}, {\em nets} and {\em webs}.
\claim{\bf Counting pencils on $\PP^1$.}
The enumerative geometry of pencils on $\PP^1$ with prescribed ramifications is ruled  by the intersection theory on $G(2,n+2)$. By~\cite{GatSant}, $A^*(G(2,n+2))$ is generated by $D_1,D_2$.
\bclm{\bf Proposition.} {\em Let $a,b\geq 0$ and $i_1,i_2\geq 1$ such that $a+2b=2n-i_1-i_2-3$. Then:
\[
\int_{n+2}D_1^aD_2^b(\ep^{i_1}\w\ep^{i_2})=\sum_{\beta=0}^b{b\choose \beta}\omega_{i_1+\beta,i_2+2b-2\beta}\label{eq:goldberg}
\]
}
\eclm
\proof In fact:
\begin{small}
\[
D_1^aD_2^b(\ep^{i_1}\w\ep^{i_2})=D_1^a\sum_{\beta=0}^b{b\choose \beta}D_1^\beta(\ep^{i_1+\beta}\w D_2^{b-\beta}\ep^{i_2})=\sum_{\beta=0}^b{b\choose \beta}D_1^{a+\beta}(\ep^{i_1+\beta}\w\ep^{i_2+2b-2\beta}).
\]
\end{small}
Taking integrals one obtains precisely formula~(\ref{eq:goldberg}).
\qed

If $i_1=1$, $i_2=2$ and $b=0$ (i.e. $a=2n$):
\[
\int _{n+2}D_1^{2n}(\ep^1\w\ep^2)=\omega_{1,2}={2n\choose n}-{2n\choose n+1}={1\over n+1}{2n\choose n},
\]
the degree of the grassmannian $G(2,n+2)$, called {\em Goldberg's formula}  in~\cite{Osserman}.
 Let now $q_1,\ldots,q_{h}$ be non negative integers not bigger than $n$,  such that  $\sum q_j=2n$.
 \bclm{{\bf Theorem}}\label{thmscherb}~\cite[Scherbak]{Scherb1}{\bf .} {\em The following formula holds:
 \be
 \int_{G(2,n+2)}\sigma_{q_1}\cdot\ldots\cdot\sigma_{q_h}\cap[G(2,n+2)]=\sum_{I\subseteq \{1,\ldots,h\}}(-1)^{h+1-|I|}{\sum_{i\in I}q_i+|I|-n-1\choose h-2}.\label{eq:forscherb}
 \ee
 }
 \eclm

 \proof
Notation as in Section~\ref{convbin}. Diagram~(\ref{eq:diag1}) guarantees the following equality:
 \[
  \int_{G(2,n+2)}\sigma_{q_1}\cdot\ldots\cdot\sigma_{q_h}\cap[G(2,n+2)]=\int_{n+2}D_{q_1}\cdot\ldots\cdot D_{q_h}(\ep^1\w\ep^2),
 \]
 where the last hand side is the coefficient of $\ep^{n+1}\w\ep^{n+2}$ in the expansion of \linebreak $D_{q_1}\cdot\ldots\cdot D_{q_h}(\ep^1\w\ep^2)$. Iterating formula~(\ref{eq:hleibrule}), one easily gets:
\[
D_{q_1}\cdot\ldots\cdot
D_{q_{h}}(\ep^1\w\ep^2)=\sum_{\tiny\matrix{1\leq j\leq h\cr 0\leq m_j\leq q_j}}\ep^{1+m_1+\cdots+m_h}\w\ep^{2+q_1+\cdots+q_h-m_1-\cdots-m_h}=
\]
 \be
 =\sum_{\tiny\matrix{1\leq j\leq h\cr 0\leq m_j\leq q_j}}\ep^{1+m_1+\cdots+m_h}\w\ep^{2+2n-m_1-\cdots-m_h}.\label{eq:sumt}
\ee
 Since $\sum_{1\leq j\leq h} q_j=2n$, the only surviving terms in the sum~(\ref{eq:sumt}) are those for which either $m_1+\cdots+m_h=n$ or  $m_1+\cdots+m_h=n+1$.
 For each $m\geq 0$, $h\geq 1$ and $1\leq j\leq h$, let
 $p_{h;j}(m):=\{(m_1,\ldots, m_h)\in p_h(m)\,|\,  m_j\leq q_j\}.
$
Then:
   \[
\int_{n+2}D_{q_1}\cdot\ldots\cdot D_{q_{{
h}}}(\ep^1\w\ep^2)=|\cap_{1\leq j\leq h}p_{h;j}(n)|-|\cap_{1\leq
j\leq h} p_{h;j}(n+1)|,
 \]
and our task consists now in evaluating the right hand side. As done in~\cite{SV}, in the context of  representation theory, we  apply the {\em inclusion-exclusion formula} (see e.g.~\cite[p.~76]{Cameron}) to our situation.  First, one notices that for each $1\leq j\leq h$, the set $p_{h;j}(m)$ is the complement in $p_h(m)$ of
$
 p'_{h;j}(m):=\{(m_1,\ldots, m_h)\in p_h(m)\,|\,m_j\geq q_j+1\}.
$
Therefore:
\[
|\cap_{1\leq j\leq h} p_{h;j}(m)|=|p_h(m)-\cup_{1\leq j\leq h} p'_{h;j}(m)| =\sum_{J\subseteq \{1,\ldots,h\}}(-1)^{|J|}|\cap_{j\in J}p'_{h;j}(m)|,
\]
where we used  De Morgan's laws and, for the second equality, the  inclusion exclusion formula, with the convention that  $|\cap_{j\in J}p'_{h;j}(m)|=|p_h(m)|$, if $|J|=\emptyset$.
Now, for each $(m_1,\ldots, m_h)\in \cap_{j\in J} \,p'_{h;j}(m)$ and each $j\in J$, replace the element $m_j$  with $m'_j=m_j-q_j-1$ to see that
$
\cap_{j\in J} \,p'_{h;j}(m)=p_h\left(m-\sum_{j\in J}q_j-|I|\right).
$
Then:
\be
|\cap_{j\in J}p'_{h;j}(m)|=\left|p_h\left(m-\sum_{j\in J}q_j-|J|\right)\right|={m-\sum_{j\in J}q_j-|J|+h-1\choose h-1},\label{eq:apines}
\ee
where to get last equality in~(\ref{eq:apines}) we applied formula~(\ref{eq:ines}).

Finally:
\[
\int_{n+2}D_{q_1}\cdot\ldots\cdot D_{q_{{
h}}}(\ep^1\w\ep^2)=|\cap_{1\leq j\leq h}p_{h;j}(n)|-|\cap_{1\leq
j\leq h} p_{h;j}(n+1)|=
\]
\[
 \sum_{J\subseteq \{1,\ldots,h\}}(-1)^{|J|}\left[{n-\sum_{j\in J}q_j-|J|+h-1\choose h-1}- {n+1-\sum_{j\in J}q_j-|J|+h-1\choose h-1}\right]=
\]
\[
=- \sum_{J\subseteq \{1,\ldots,h\}}(-1)^{|J|}{n-\sum_{j\in J}q_j-|J|+h-1\choose h-2}=
\]
\[
=- \sum_{\tiny\matrix {I\subseteq
\{1,\ldots,h\}}}(-1)^{h-|I|}{\sum_{i\in I}q_i+|I|-n-1\choose h-2},
\]
%{\ros * In mia oppinione nel'ultima uguaglianza dovrebbe essere
%messo soltanto $I\subset \{1,\ldots h\}$ perche questo é
%equivalente a scrivere $\matrix{J\subseteq \{1,\ldots,h\}\cr
%I\subseteq \{1,\ldots,h\}\setminus J}$. Quando si fa la
%sostituzione $\sum_{i\in I}q_i=2n-\sum_{j\in J} q_j$ il $J$ non
%\`{e} piu rilevante. Quindi secondo me non c'e bisogno piu
%caricarlo nella formula.}
where, for $I=\{1,\ldots,h\}\setminus J$, we have used the
equality $\sum_{i\in I}q_i=2n-\sum_{j\in J} q_j$. The last term of
the equalities above obviously coincides with the r.h.s.
of~(\ref{eq:forscherb}). \qed \claim{\bf Remark.}
Expression~(\ref{eq:forscherb}) is Scherbak's formula in~\cite{Scherb1}, written in a { slightly} modified
version  to (formally)
include the cases $h=1$ (the degree $(=0)$ of a special Schubert cycle) and $h=2$ (the degree of the product of two special Schubert cycles).

\claim{\bf Counting nets on $\PP^1$.} When $r=2$, let $C$ be the image of $\PP^1$ through the rational map $\phi:\PP^1\sra \PP(V)$ induced the given $g^2_d$. Assume that $\phi$ is a morphism, i.e. the $g^2_d$ has no base point.  The geometrical interpretation of the $V$-ramification sequence at a point $P$ in terms of the nature of the point $Q:=\phi(P)$ on $C$ is as follows. At a general point  of $\PP^1$, the $V$-ramification sequence   is $(1,2,3)$ (weight $0$) and $Q$ is an ordinary point of $C$. Instead, $Q$  is a {\em flex}, a {\em  hyperflex}, a {\em  cusp}, or a {\em   tacnode} if the $V$-ramification sequence at $P$ is respectively $(1,2,4)$, $(1,2,5)$, $(1, 3,4)$, $(1,3,5)$. The $V$-weights are respectively $1$, $2$, $2$ and $3$. Suppose that $a,b,c,d$ are four integers such that $a+2b+2c+3d=3n$. Then by~\ref{recalEH},~(\ref{eq:giamb}) and~(\ref{eq:diag1}), with the same notation as in Section~\ref{recallsc}:
\begin{eqnarray*}
N_{a,b,c,d}&=&\int_{G(3,n+3)}\sigma_1^a\cdot \sigma_2^{b}\cdot\sigma_{(134)}^c\cdot\sigma_{(135)}^d\cap[G(3,n+3)]=\\
&=&\int_{n+3}D_1^a\cdot D_2^{b}\cdot(\Delta_{(134)}(D))^c\cdot(\Delta_{(135)}(D))^d(\ep^1\w\ep^2\w\ep^3)
\end{eqnarray*}
is the number  of  projectively non equivalent rational plane curves of degree $n+2$ having $a$ flexes, $b$ hyperflexes, $c$ cusps and $d$ tacnodes at $a+b+c+d$  distinct points. By definition of $\Delta$ (Section~\ref{gradzt}) and of $D_t^{-1}$ (Section~\ref{clm28}):
\be
\Delta_{(134)}(D)=D_1^2-D_2=\ovD_2\quad \mathrm{and}\quad  \Delta_{(135)}(D)=D_1D_2-D_3=D_1\ovD_2-\ovD_3\label{eq:d3d3bar}
\ee
In fact, $D_t\circ D_t^{-1}=1$ implies that $D_3=D_2\ovD_1-D_1\ovD_2+\ovD_3$ and since $\ovD_1=D_1$ one gets last member of the second equality of~(\ref{eq:d3d3bar}).
Thus:
\[
N_{a,b,c,d}=\int_{n+3}D_1^a\cdot D_2^{b}\cdot\ovD_2^c\cdot(D_1\ovD_2-\ovD_3)^d(\ep^1\w\ep^2\w\ep^3)=
\]
\[
=\sum_{d'=0}^d(-1)^{d'}\int_{n+3}{d\choose d'}D_1^{a+d'}D_2^b\,\ovD_2^{c+d'}\,\ovD_3^{d-d'}(\ep^1\w\ep^2\w\ep^3)
\]
\bclm{\bf Theorem.} The number $N_{a,b,c,d}$ is given by:
\[
N_{a,b,c,d}=\sum_{\tiny\matrix{0\leq d'\leq d\cr{\bm\beta}\in p_4(b)\cr {\bm\gamma}\in p_3(c+d')}}(-1)^{d'}{d\choose d'}{b\choose {\bm\beta}}{ c+d'\choose {\bm\gamma}}\,\omega_{(i_1(d,{\bm\beta},{\bm\gamma}),i_2(d,{\bm\beta},{\bm\gamma}),i_3(d,{\bm\beta},{\bm\gamma}))}
\]
where
\be
\left\{\matrix{i_1(d,{\bm\beta},{\bm\gamma})&=&1+d-d'+\gamma_2+\gamma_3+\beta_1\cr
i_2(d,{\bm\beta},{\bm\gamma})&=&2+d-d'+\gamma_1+\gamma_3+\beta_2+2\beta_4\cr
i_3(d,{\bm\beta},{\bm\gamma})&=&3+d-d'+\gamma_1+\gamma_2+2\beta_3+\beta_2}\right.\label{eq:indplanecurv}
\ee
\eclm
\proof
By applying the formulas proven in  Section~\ref{SchubCalc}. First one computes:
\[
\ovD_2^{c+d'}\ovD_3^{d-d'}(\ep^1\w\ep^2\w\ep^3)=\ovD_2^{c+d'}(\ep^{1+d-d'}\w\ep^{2+d-d'}\w\ep^{3+d-d'})
\]
\[
=\sum_{{\bm\gamma}\in p_3(c+d')}{c+d'\choose {\bm\gamma}}\ep^{1+d-d'+\gamma_2+\gamma_3}\w\ep^{2+d-d'+\gamma_1+\gamma_3}\w\ep^{3+d-d'+\gamma_1+\gamma_2},
\]
where the first equality is by~(\ref{eq:iterdbar}) and the second one is by Proposition~\ref{dbarhm1}.
Moreover
\[
D_1^{a+d}D_2^b(\ep^{1+d-d'+\gamma_2+\gamma_3}\w\ep^{2+d-d'+\gamma_1+\gamma_3}\w\ep^{3+d-d'+\gamma_1+\gamma_2})=
\]
\[
=\sum {b\choose {\bm\beta}} D_1^{a+d+\beta_1}(\ep^{i_1(d,{\bm\beta},{\bm\gamma})}\w \ep^{i_2(d,{\bm\beta},{\bm\gamma})}\w \ep^{i_3(d,{\bm\beta},{\bm\gamma}))})=\sum {b\choose {\bm\beta}}\omega_{(i_1(d,{\bm\beta},{\bm\gamma}),i_2(d,{\bm\beta},{\bm\gamma}),i_3(d,{\bm\beta},{\bm\gamma}))}
\]
with $i_1(d,{\bm\beta},{\bm\gamma})$, $i_2(d,{\bm\beta},{\bm\gamma})$, $i_3(d,{\bm\beta},{\bm\gamma})$ given by (\ref{eq:indplanecurv}). Putting all together, one gets precisely the claimed formula for $N_{a,b,c,d}$
\qed

\claim{\bf Counting webs on $\PP^1$.} \label{countingwebs} Let $V$
be a base point free  web on $\PP^1$ and let $\phi:\PP^1\sra
\PP(V)$ be the induced morphism. Let $P\in \PP^1$. We shall say
that $Q:=\phi(P)$ is a {\em stall}, a {\em hyperstall}, a {\em
flex}, a {\em cusp}, if the $V$-ramification sequence at $P$ is,
respectively, $(1,2,3,5)$, $(1,2,3,6)$, $(1,2,4,5)$, $(1,3,4,5)$.
At a {\em stall} (resp. at {\em a hyperstall}) $Q\in C$ the
osculating plane meets the curve at $Q$ with multiplicity $4$
(resp. with multiplicity $5$). If $Q$ is a flex, all the planes of
the pencil containing the tangent line to $Q$ meets the curve with
multiplicity at least $3$. If $Q$ is a cusp, the local {
analytic} equation of $C$ at $Q$ is precisely $y^2-x^3=0$. Imitating
what we did for rational plane curves, we may write a formula
counting the number (with multiplicity) of all the (projectively
non equivalent) rational space curves having $a$ stalls, $b$
hyperstalls, $c$ flexes and $d$ cusps, at $a+b+c+d$ distinct
points  such that $a+2b+2c+3d=4n$. Such a
number is counted by the integral (Cf.~\ref{recalEH},
~\ref{recallsc},~\ref{clm28}):
\[
f_{a,b,c,d}:=\int_{G(4,n+4)}\sigma_1^a\sigma_2^b\sigma_{1245}^c\sigma_{1345}^d\cap [G(4,n+4)]=\int_{n+4} D_1^aD_2^b\ovD_2^c\ovD_3^d(\ep^1\w\ep^2\w\ep^3\w\ep^4)
\]
i.e. the coefficient of  $\pi_{4,n+4}$ (see Section~\ref{secpkn}),  in the expansion of
\be
D_1^aD_2^b\ovD_2^c\ovD_3^d(\ep^1\w\ep^2\w\ep^3\w\ep^4)
\label{eq:prenabcd}
\ee
\bclm{\bf Theorem.} {\em The number $f_{a,b,c,d}$ is given by:
\be
f_{a,b,c,d}=\sum_{\tiny\matrix{{\bm\beta}\in p_5(b)\cr {\bm \gamma}\in p_4(c)\cr {\bm\delta}\in p_4(d)\cr 0\leq l\leq \beta_2+\beta_3\cr 0\leq m\leq \gamma_1+\gamma_2}}{b\choose{\bm\beta}} {c\choose {\bm\gamma}}{d\choose{\bm\delta}}{\beta_2+\beta_3\choose l}{\gamma_1+\gamma_2\choose m}\omega_{I({\bm\beta},{\bm\gamma},{\bm\delta};l,m)}\label{eq:fabcd}
\ee
having set:
\be
I({\bm \beta},{\bm \gamma},{\bm\delta};l,m)=(i_1({\bm \beta},{\bm \gamma},{\bm\delta};l,m),i_2({\bm \beta},{\bm \gamma},{\bm\delta};l,m), i_3({\bm \beta},{\bm \gamma},{\bm\delta};l,m),i_4({\bm \beta},{\bm \gamma},{\bm\delta};l,m))\label{eq:ibgdlm}
\ee
with
 \be
\left\{{\small\matrix{i_1(\beta,\gamma,\delta;l,m)&=&1\,+&\sum_{j\neq 1}\delta_j+\gamma_1+\gamma_4+\beta_1;\hskip166pt\cr{}\cr
i_2(\beta,\gamma,\delta;l,m)&=&2\,+&\sum_{j\neq 2}\delta_j+\gamma_2+\gamma_4+\beta_2+2\beta_5;\hskip130pt\cr{}\cr
i_3(\beta,\gamma,\delta;l,m)&=&3\,+&\sum_{j\neq 3}\delta_j+\gamma_3+\beta_3+l+m;\hskip160pt\cr{}\cr
i_4(\beta,\gamma,\delta;l,m)&=&4\,+&\sum_{j\neq 4}\delta_j+\gamma_1+\gamma_2+\gamma_3+\beta_2+2\beta_4+\beta_3-l-m\hskip37pt}
}\right.
\ee

\medskip
\noindent
(where $\omega_{(\beta,\gamma,\delta;l,m)}=0$ if $i_p(\beta,\gamma,\delta;l,m)>n+4$ for  some $p\in\{1,2,3,4\}$).
}
\eclm
\proof
The proof is straightforward. By formula~(\ref{eq:lemcomp}), on first has:
\[
D_1^aD_2^b\ovD_2^c\ovD_3^d(\ep^1\w\ep^2\w\ep^3\w\ep^4)=
\]
\[
=D_1^aD_2^b\ovD_2^c\sum_{\delta\in p_4(d)}{d!\over \delta!}\ep^{1+\sum_{j\neq 1}\delta_j}\w \ep^{2+\sum_{j\neq 2}\delta_j}\w \ep^{3+\sum_{j\neq 3}\delta_j}\w \ep^{4+\sum_{j\neq 4}\delta_j}=
\]
One now applies~(\ref{eq:forovD2m}) getting:
\[
=D_1^aD_2^b\sum_{\tiny\matrix{\delta\in p_4(d)\cr{\bm \gamma}\in p_4(c)\cr 0\leq m\leq \gamma_1+\gamma_2}}{c\choose {\bm\gamma}}{d\choose{\bm\delta}}{\gamma_1+\gamma_2\choose m}\ep^{i_1({\bm \gamma},{\bm\delta};m)}\w\ep^{i_2({\bm \gamma},{\bm\delta};m)}\w\ep^{i_3({\bm\gamma,\delta};m)}\w\ep^{i_4({\bm \gamma},{\bm\delta};m)},
\]
where
\[
\left\{\matrix{i_1({\bm \gamma},{\bm\delta};m)&=&1&+\sum_{j\neq 1}\delta_j+\gamma_1+\gamma_4\cr
i_2({\bm \gamma},{\bm\delta};m)&=&2&+\sum_{j\neq 2}\delta_j+\gamma_2+\gamma_4\cr
i_3({\bm \gamma},{\bm\delta};m)&=&3&+\sum_{j\neq 3}\delta_j+ \gamma_3+m\cr
i_4({\bm \gamma},{\bm\delta};m)&=&4&+\sum_{j\neq 4}\delta_j+\gamma_1+\gamma_2+\gamma_3-m
}\right.
\]
Finally, applying formula~(\ref{eq:forD2m}):
\[
=\sum_{\tiny\matrix{{\bm\beta}\in p_5(b)\cr {\bm\gamma}\in p_4(c)\cr {\bm\delta}\in p_4(d)\cr 0\leq l\leq \beta_2+\beta_3\cr 0\leq m\leq \gamma_1+\gamma_2}}{b\choose {\bm\beta}} {c\choose{\bm\gamma}} {d\choose {\bm\delta}}{\beta_2+\beta_3\choose l}{\gamma_1+\gamma_2\choose m} D_1^{a+\beta_1}\bfep^{I(\beta,\gamma,\delta;l,m)},
\]
where $I({\bm\beta},{\bm \gamma},{\bm\delta};l,m)$ is given precisely by~(\ref{eq:ibgdlm}). Taking integrals, one gets precisely formula~(\ref{eq:fabcd}).\qed

\claim{}\label{finalsect} Let $a=c=d=0$. Then $b=2n$. For each ${\bm\beta}\in p_5(2n)$, let:
\[
I({\bm\beta},l)=(1+\beta_1 ,  2+\beta_2+2\beta_5,
3+\beta_3+l, 4+\beta_2+\beta_3+2\beta_4-l).
\] Then:
\[
HS_n:=f_{0,2n,0,0}=\sum_{\tiny\matrix{\beta\in p_5(2n)\cr 0\leq l\leq \beta_2+\beta_3}} {2n\choose {\bm \beta}}{\beta_2+\beta_3\choose l}\omega_{I({\bm\beta};l)},
\]
is the number (with multiplicities) of rational space curves having $2n$ hyperstalls at $2n$ prescribed distinct points.

Similarly, if $a,b,d=0$, then $c=2n$. For each ${\bm\gamma}\in p_4(2n)$, let
\[
I({\bm\gamma};m)=(1+\gamma_1+\gamma_4, \,2+\gamma_2+\gamma_4, \,3+\gamma_3+m, \,4+\gamma_1+\gamma_2+\gamma_3-m)
\]
Thus:
\[
f_{0,0,2n,0}=\sum_{\tiny\matrix{{\bm\gamma}\in p_4(2n)\cr 0\leq m\leq \gamma_1+\gamma_2}} {2n\choose {\bm \gamma}}{\gamma_1+\gamma_2\choose m}\omega_{I({\bm \gamma};m)},
\]
which is our way to answer to Ranestad's original question.

\bigskip
\noindent

%\claim{} (See~\cite{donagi}, p.~799--800, \cite{GH}) The number of lines in $\PP^{2n+1}$ intersecting four given $n$-spaces in general position is computed by
%\[
%L_n:=\int_{G(2,2n+2)}\sigma_n^4\cap[G(2,2n+2)]
%\]
%Applying our dictionary, this is the same as computing the coefficient of $\ep^{1+2n}\w\ep^{2+2n}$ in the expansion of $D_n^4(\ep^1\w\ep^2)\in\w^2M_{2+2n}$.
%One has:
%\[
% L_n:=D_n^4(\ep^1\w\ep^2)=D_n^2( D_n^2(\ep^1\w\ep^2))=D_n^2(\ep^1\w\ep^{2+2n}+ 2D_{n-1}(\ep^2\w\ep^{2+n})+ D_{n-1}^2(\ep^3\w\ep^2))
% \]
%where in the last equality we used the Newton formula for $h=n$ and $m=2$. Applying~(\ref{eq:preip}) and the skew-symmetry of $\w$, last member can be written as:
% \[
% D_n^2(\ep^{1}\w\ep^{2+2n}+2D^2_{n-1}(\ep^{2}\w\ep^{3})- D_{n-1}^2(\ep^2\w\ep^3))= D_n^2(\ep^{1}\w\ep^{2+2n}+D_{n-1}^2(\ep^2\w\ep^3))=
%\]
%i.e.
%\[
%=D_n^2(\ep^{1}\w\ep^{2+2n})+D_n^2D_{n-1}^2(\ep^2\w\ep^3)=\pi_{2,2+2n}+D_{n-1}^2D_n^2(\ep^2\w\ep^3)
%\]
%By applying again Newton's formula, one gets: $D_n^2(\ep^2\w\ep^3)=D_{n-1}^2(\ep^3\w\ep^4)$, so that:
%\[
%L_n=\pi_{2,2+2n}+ D_{n-1}^4(\ep^3\w\ep^4)=(1+L_{n-1})\pi_{2,2+2n}.
%\]
%Therefore:
%$
%L_n=1+L_{n-1}
%$. Since $L_0=1$, it follows that $L_n=1+n$.

\bigskip
\noindent
Dipartimento di Matematica, Politecnico di Torino, C.so Duca degli Abruzzi 24, 10129 (TO), Italy

\smallskip
\noindent
Departamento de Ci\^{e}ncias Exatas, U. E. de Feira de Santana, (BA) - 44031-460
Brazil

\medskip
{
\noindent {\tt cordovez@calvino.polito.it}\\
\noindent {\tt letterio.gatto@polito.it}\\
\noindent {\tt taisesantiago@uefs.br} }
\end{document}